\documentclass[oneside,american]{amsart}
\usepackage{mathptmx}

\usepackage[T1]{fontenc}
\usepackage[latin9]{inputenc}
\usepackage{geometry}
\usepackage{xcolor}
\usepackage{verbatim}
\usepackage{amsbsy}
\usepackage{amstext}
\usepackage{amsthm}
\usepackage{amssymb}
\usepackage{graphicx}
\PassOptionsToPackage{normalem}{ulem}
\usepackage{ulem}

\makeatletter

\providecommand{\tabularnewline}{\\}
\providecolor{lyxadded}{rgb}{0,0,1}
\providecolor{lyxdeleted}{rgb}{1,0,0}
\DeclareRobustCommand{\lyxsout}[1]{\ifx\\#1\else\sout{#1}\fi}

\numberwithin{equation}{section}
\numberwithin{figure}{section}
\theoremstyle{plain}
\newtheorem{thm}{\protect\theoremname}
\theoremstyle{remark}
\newtheorem{rem}[thm]{\protect\remarkname}

\makeatother

\usepackage{babel}
\providecommand{\remarkname}{Remark}
\providecommand{\theoremname}{Theorem}

\begin{document}
\title{Adaptive algorithm for electronic structure calculations using reduction
of Gaussian mixtures}
\author{Gregory Beylkin, Lucas Monz\'{o}n and Xinshuo Yang }
\address{Department of Applied Mathematics \\
 University of Colorado at Boulder \\
 UCB 526 \\
 Boulder, CO 80309-0526}
\begin{abstract}
We present a new adaptive method for electronic structure calculations
based on novel fast algorithms for reduction of multivariate mixtures.
In our calculations, spatial orbitals are maintained as Gaussian mixtures
whose terms are selected in the process of solving equations.

Using a fixed basis leads to the so-called ``basis error'' since
orbitals may not lie entirely within the linear span of the basis.
To avoid such an error, multiresolution bases are used in adaptive
algorithms so that basis functions are selected from a fixed collection
of functions, large enough to approximate solutions within any user-selected
accuracy.

Our new method achieves adaptivity without using a multiresolution
basis. Instead, as part of an iteration to solve nonlinear equations,
our algorithm selects the ``best'' subset of linearly independent
terms of a Gaussian mixture from a collection that is much larger
than any possible basis since the locations and shapes of the Gaussian
terms are not fixed in advance. Approximating an orbital within a
given accuracy, our algorithm yields significantly fewer terms than
methods using multiresolution bases.

We demonstrate our approach by solving the Hartree-Fock equations
for two diatomic molecules, $\mbox{HeH}^{+}$ and $\mbox{LiH}$, matching
the accuracy previously obtained using multiwavelet bases.

.
\end{abstract}

\maketitle

\section{\label{sec:Introduction}Introduction}

We present a new adaptive method for electronic structure calculations
based on algorithms for reduction of multiresolution multivariate
mixtures \cite{BE-MO-YA:2019}. While we represent solutions using
a linear combination of Gaussians (which has a long history in quantum
chemistry), these Gaussians are not selected in advance to be used
as an approximate basis but are generated in the process of solving
equations. Thus, we avoid the so-called ``basis error'' usually
associated with methods using Gaussians. Our approach can also be
characterized as a ``gridless'' or a ``meshfree'' method.

Using Gaussians to find solutions of quantum chemistry problems has
its origins in seminal papers \cite{BOYS:1960,LON-SIN:1960,SINGER:1960}
and is motivated by the fact that integrals involving these functions
can be evaluated efficiently. In these early quantum chemistry papers
the authors proposed to use linear combinations of Gaussians whose
exponents and coefficients were found (or optimized) via Newton's
method in order to capture the correct behavior near the nuclear cusps
and the correct rate of decay. However, this approach proved unsustainable
as problems became larger. Instead, the construction of a basis for
spatial orbitals has been performed off-line and the resulting sets
of functions were then used as a fixed basis, leading to the so-called
``basis error'' if the actual solution is not well approximated
within the linear span of such a basis.

The use of these Gaussian bases have revolutionized computational
quantum chemistry in spite of the absence of a systematic way for
controlling accuracy or providing guaranteed error bounds. In fact
selecting an appropriate basis set became an art form requiring insight
into the underlying solution. However, once a basis set is selected,
the accuracy of the solution obtained using such a basis is ultimately
limited. Since the equations being solved, e.g. the Hartree-Fock equations
or the Kohn\textendash Sham equations of density functional theory
(DFT) (see e.g. \cite{PAR-YAN:1989}), only provide approximate solutions,
the limitation in accuracy makes it difficult to separate the impact
of using approximate equations from the approximate methods for solving
them. For this reason, the use of adaptive methods to avoid the loss
of accuracy caused by basis sets is highly desirable in this field.

The advent of multiresolution analysis (see e.g. \cite{DAUBEC:1992,BE-CO-RO:1991})
laid a conceptual foundation for adaptive methods but it took some
time before practical adaptive algorithms were developed using multiwavelets
\cite{ALPERT:1993}; see \cite{H-F-Y-G-B:2004,Y-F-G-H-B:2004,Y-F-G-H-B:2004a}.
A central element for the success of the multiresolution algorithms
can be traced to the fact that physically significant (integral) operators
arising in problems of quantum chemistry are naturally represented
by radial kernels, which, in turn, can be accurately approximated
by a linear combination of Gaussians. The key advantage of using such
approximations is that they yield a separated representation in which
operators are efficiently applied one direction at a time. Without
a separated representation, multiresolution operators would be too
expensive to be practical in dimensions three and higher.

Multiresolution methods systematically refine basis functions (i.e.,
numerical grids) in the vicinity of the cusp-type singularities while
using relatively few basis functions elsewhere. These methods have
proven successful in efficiently computing highly accurate solutions,
achieving guaranteed error bounds and eliminating the basis error.
This approach has been implemented at ORNL in the software package
MADNESS ( Multiresolution ADaptive Numerical Environment for Scientific
Simulation, see \cite{H-B-B-C-F-F-G-etc:2016}), which is now considered
the most accurate approach in this field \cite{J-S-F-H-B-G-2017}.
Wavelets are also used in BigDFT \cite{G-N-G-D-G-W-C-Z-R-B:2008}
(cf. a mixed-basis method using plane waves and atom-centered radial
polynomials \cite{PAH-HAN:2002} and interlocking multi-center grids
\cite{SHI-HIR:2007}, as examples).

However, since wavelets (or multiwavelets) do not resemble the spatial
orbitals, a large number of basis functions is needed to represent
solutions, e.g. an individual orbital may require $\approx2\cdot10^{5}$
basis functions in three dimensions. Moreover, such local refinement
schemes do not take advantage of the essential simplicity of the spatial
orbitals far from the nuclei and require boundary conditions to limit
the computational domain. While adaptive multiresolution methods are
sufficiently fast to be used within one-particle theories of quantum
chemistry, an advancement towards two particle theories or solving
Schrödinger's equation had a limited success due to computational
costs (see e.g. \cite{BE-MO-PE:2007,BE-MO-PE:2008}). 

Hence it is of interest to develop new adaptive schemes and compare
them with MADNESS. As it was demonstrated in \cite{BEY-HAU:2013},
it is possible to iteratively solve equations of quantum chemistry
using new algorithms based on nonlinear approximation of functions
to reduce the number of terms in intermediate representations (without
resorting to Newton's method). The results of this approach are presented
in \cite{BEY-HAU:2013} and are based on using Slater-type orbitals. 

In this paper we present a new adaptive method that uses linear combinations
of Gaussians to represent the solutions. We demonstrate the new approach
by solving the Hartree-Fock equations for $\mbox{HeH}^{+}$ and $\mbox{LiH}$
so that we can compare the resulting representations with those in
\cite{H-F-Y-G-B:2004,BEY-HAU:2013}. As in \cite{H-F-Y-G-B:2004,BEY-HAU:2013},
we formulate the problem using integral equations and use a convergent
iteration to solve them. In contrast to the above approaches, we use
linear combinations of Gaussians to accurately approximate not only
operators and potentials but also the functions on which they operate.
As a result, all integrals can be computed explicitly and exactly
by simply updating the parameters of the Gaussians involved. The computational
effort thus moves from that of approximating and computing integrals
to that of maintaining a reasonable number of terms in intermediate
representations of solutions within the iterative scheme. Using our
approach, computing a single integral (e.g. involving a Green's function,
potential and a wave function) can easily yield a Gaussian mixture
with $\approx10^{6}$ terms, most of which are linearly dependent
(within a user-selected accuracy). We use the reduction algorithm
from \cite[Algorithm 1]{BE-MO-YA:2019} to reduce the number of terms
after each operation by finding the ``best'' linearly independent
terms, thus maintaining a reasonably small number of them. This algorithm
has computational complexity $\mathcal{O}\left(r^{2}N+p\left(d\right)rN\right)$,
where $N$ is the original number of terms and $r$ is the number
of so-called skeleton terms selected by the algorithm and $p\left(d\right)$
is the cost of computing the inner product between Gaussians as a
function of dimension (in this paper $d=3$ and $p\left(d\right)$
is a small constant). There is no underlying grid to maintain (thus,
``gridless'' or ``meshfree'' method) and there is no need to impose
boundary conditions to reduce the size of the computational domain
as in \cite{H-F-Y-G-B:2004}. 

Since we only present results using examples of two small molecules,
it is natural to ask if our approach generalizes to large molecules.
The obvious obstacle is the quadratic dependence of the reduction
algorithm on the number of skeleton terms (clearly, the final number
of terms to represent solutions for large molecules is expected to
be large). The key to addressing this problem is to avoid the unnecessary
application of the reduction algorithm to sets of functions that are
obviously linearly independent. Algorithmically, it means splitting
terms of a Gaussian mixture into groups that are likely to have linear
dependence (and never perform a global reduction). This approach allows
us to maintain a reasonable number of skeleton terms within each group.
Moreover, the resulting reduction algorithm is trivially parallel
since the reduction within each group is independent of other groups.
We demonstrate subdivision into groups using both location and scale
in our examples and show that the accuracy of the solution improves
gradually due to a convergent iteration.

We also observe that the Hartree\textendash Fock setup does not take
into account the electron-electron cusps whereas it is well known
that incorporating them significantly improves the accuracy of energy
calculations. In our approach, it would be natural to introduce a
linear combination of Gaussian geminals (originally proposed in the
seminal papers \cite{BOYS:1960,SINGER:1960}). Using Gaussian geminals
to improve accuracy has a long history in quantum chemistry (see \cite{KOM-KIN:2011,BAR-LOO:2017}
on this topic and references therein). We note that incorporating
geminal functional forms into the adaptive approach of this paper
appears possible (but by no means trivial) since the reduction algorithm
for multivariate mixtures discussed in \cite{BE-MO-YA:2019} (and
on which the approach of this paper is based) is applicable to much
more general functions than those used in our two examples. We also
note that a different approach to incorporate cusps and to address
basis set error was recently proposed in \cite{G-P-F-A-S-T:2018}
by using a range-separated DFT formulation. We plan to address adaptivity
involving electron-electron interactions and to demonstrate our approach
on large molecules elsewhere.

We start by briefly stating the Hartree\textendash Fock equations
in Section~\ref{sec:Equations-of-quantum-qc} and demonstrating that
the functional form of our approximate solution, a Gaussian mixture,
is maintained when solving the integral form of these equations via
iteration (we note that the same reasoning applies to the Kohn-Sham
equations). Then we briefly review approximations of operators and
potentials that we use in Section~\ref{sec:Examples}. We then turn
to the reduction algorithm in Section~\ref{sec:Reduction-algorithm}
and describe a subdivision scheme which is critical to the practical
use of our approach. We present examples and comparisons in Section~\ref{sec:Examples}
and discuss the results in Section~\ref{sec:Conclusions-and-further-work}.

\section{\label{sec:Equations-of-quantum-qc}Solving the Hartree-Fock equations }

We use an integral version of the Hartree\textendash Fock equations
(see e.g. \cite{SZA-OST:1996}) and present these equations for the
reader's convenience in order to emphasize a combination of properties
that make our approach possible. We note that our approach is equally
applicable to the Kohn-Sham equations considered in \cite{H-F-Y-G-B:2004}.
Our method is based on computing orbitals as Gaussian mixtures and
using the reduction algorithm in \cite[Algorithm 1]{BE-MO-YA:2019}
to keep the number of terms under control. 

Briefly, the occupied orbitals are the lowest $N$ eigenfunctions
of the Hartree\textendash Fock operator 
\begin{equation}
\mathcal{F}\phi_{j}\left(\mathbf{r}\right)=E_{j}\phi_{j}\left(\mathbf{r}\right),\ \ j=1,\dots N,\label{eq: Fock operator}
\end{equation}
where $\mathcal{F}=-\frac{1}{2}\Delta+\mathcal{V}_{tot}$, $\mathcal{V}_{tot}=\mathcal{V}_{ext}+2\mathcal{J}-\mathcal{K}$.
The eigenvalues $E_{j}<0$ are referred to as orbital energies and
are negative for the occupied orbitals. The external potential $\mathcal{V}_{ext}$
accounts for attraction of electrons to the nuclei at locations $\mathbf{R}_{l}$,
$l=1,\dots L$,
\begin{equation}
\left(\mathcal{V}_{ext}\phi\right)\left(\mathbf{r}\right)=\left(\sum_{l=1}^{L}\frac{Z_{l}}{\left\Vert \mathbf{r}-\mathbf{R}_{l}\right\Vert }\right)\phi\left(\mathbf{r}\right),\,\,\,Z_{m}<0,\label{eq: external potential}
\end{equation}
where $L$ is the number of nuclei, $Z_{m}$ is the charge of the
$m$th nucleus and $\left\Vert \cdot\right\Vert $ is the standard
Euclidean vector norm. The Coulomb operator $\mathcal{J}$ describes
the potential created by orbitals $\phi_{i}\left(\mathbf{r}\right)$,
\[
\left(\mathcal{J}\phi\right)\left(\mathbf{r}\right)=\phi\left(\mathbf{r}\right)\left(-4\pi\Delta^{-1}\left(\sum_{i=1}^{N}\left|\phi_{i}\left(\mathbf{r}\right)\right|^{2}\right)\right),
\]
and the exchange operator $\mathcal{K}$ is defined as 
\[
\left(\mathcal{K}\phi\right)\left(\mathbf{r}\right)=\sum_{i=1}^{N}\phi_{i}\left(\mathbf{r}\right)\left(-4\pi\Delta^{-1}\left(\phi_{i}^{*}\phi\right)\right),
\]
where $\phi^{*}$ indicates the complex conjugate of the function
$\phi$. The orbitals are obtained solving the $N$ coupled equations
(\ref{eq: Fock operator}) and the total energy, $E_{tot}$, is calculated
by adding the nucleus\textendash nucleus repulsion energies to the
total electron energy

\[
E_{tot}=\sum_{j=1}^{N}\left(E_{j}+\left\langle \left(-\frac{1}{2}\Delta+\mathcal{V}_{ext}\right)\phi_{j},\phi_{j}\right\rangle \right)+\sum_{k=1}^{L}\sum_{l>k}^{L}\frac{Z_{l}Z_{k}}{\left\Vert \boldsymbol{R}_{l}-\boldsymbol{R}_{k}\right\Vert }.
\]
 In order to obtain integral equations for (\ref{eq: Fock operator}),
we follow \cite{KALOS:1962,KALOS:1963,H-F-Y-G-B:2004,BEY-HAU:2013}
and use that the operator $-\Delta+\mu^{2}$ has Green's function
\[
\mathcal{G}_{\mu}\left(\mathbf{r}-\mathbf{r}'\right)=\frac{1}{4\pi}\frac{e^{-\mu\left\Vert \mathbf{r}-\mathbf{r}'\right\Vert }}{\left\Vert \mathbf{r}-\mathbf{r}'\right\Vert }.
\]
 That is,
\[
\left(-\Delta+\mu^{2}\right)\mathcal{G}_{\mu}\left(\mathbf{r}-\mathbf{r}'\right)=\delta\left(\mathbf{r}-\mathbf{r}'\right).
\]
Note that the kernel of $\Delta^{-1}$ in the equations above is 
\[
\mathcal{G}_{0}\left(\mathbf{r}-\mathbf{r}'\right)=\frac{1}{4\pi}\frac{1}{\left\Vert \mathbf{r}-\mathbf{r}'\right\Vert }
\]
and that (\ref{eq: Fock operator}) is equivalent to 
\begin{equation}
\phi_{j}\left(\mathbf{r}\right)=-2\mathcal{G}_{\mu_{j}}*\left(\mathcal{V}_{tot}\phi_{j}\right)\left(\mathbf{r}\right),j=1,\dots,N,\label{eq:integral H-F}
\end{equation}
where $\mu_{j}^{2}=-2E_{j}$. We solve (\ref{eq:integral H-F}) via
an iteration in which $\mu_{j}$ are changing and approach $\sqrt{-2E_{j}}$.
There are four observations (properties) that make our approach possible
as follows.
\begin{enumerate}
\item The iteration to solve (\ref{eq:integral H-F}) is convergent (see
additional comments below).
\item All potentials and Green's functions can be accurately and efficiently
approximated by linear combinations of Gaussians. This property of
potentials and Green's functions is the foundation of the approach
in MADNESS, see e.g. \cite{H-F-Y-G-B:2004,H-B-B-C-F-F-G-etc:2016}. 
\item Since in our approach we represent orbitals $\phi_{j}$ via linear
combinations of Gaussians and potentials and Green's functions are
approximated by linear combinations of Gaussians, all integrals are
evaluated explicitly.
\item In order to keep the number of terms in the representations of the
orbitals $\phi_{j}$ under control, the key tool in our approach is
a reduction algorithm that selects a subset of the ``best'' linearly
independent terms from a large number of terms which result from applying
operators. 
\end{enumerate}
The integral form of the Kohn-Sham equations has the same properties
(cf. \cite{H-F-Y-G-B:2004}, where properties 1 and 2 are essential)
so that our approach can be used for these equations as well.

\subsection{\label{subsec:Iteration-for-solving}An iteration to solve integral
equations}

We solve (\ref{eq:integral H-F}) via the following iteration. We
initialize using a collection of Gaussians centered inside the convex
polyhedron defined by nuclear centers (see Remark~\ref{rem:By-examining-operators}).

Then, at step $m$, we first apply the operator $\mathcal{V}_{tot}^{\left(m\right)}$
to the current approximation of the orbitals and compute
\begin{eqnarray}
\begin{pmatrix}\mathcal{V}_{tot}^{\left(m\right)}\phi_{1}^{\left(m\right)}\\
\mathcal{V}_{tot}^{\left(m\right)}\phi_{2}^{\left(m\right)}\\
\vdots\\
\mathcal{V}_{tot}^{\left(m\right)}\phi_{N}^{\left(m\right)}
\end{pmatrix} & \leftarrow & \begin{pmatrix}\left(\mathcal{V}_{ext}+2\mathcal{J}^{\left(m\right)}-\mathcal{K}^{\left(m\right)}\right)\phi_{1}^{\left(m\right)}\\
\left(\mathcal{V}_{ext}+2\mathcal{J}^{\left(m\right)}-\mathcal{K}^{\left(m\right)}\right)\phi_{2}^{\left(m\right)}\\
\vdots\\
\left(\mathcal{V}_{ext}+2\mathcal{J}^{\left(m\right)}-\mathcal{K}^{\left(m\right)}\right)\phi_{N}^{\left(m\right)}
\end{pmatrix}\label{eq:HF-iteration}
\end{eqnarray}
where 
\[
\left(\mathcal{V}_{ext}\phi_{j}^{\left(m\right)}\right)\left(\mathbf{r}\right)=\left(\sum_{l=1}^{L}\frac{Z_{l}}{\left\Vert \mathbf{r}-\mathbf{R}_{l}\right\Vert }\right)\phi_{j}^{\left(m\right)}\left(\mathbf{r}\right),
\]
\[
\left(\mathcal{J}^{\left(m\right)}\phi_{j}^{\left(m\right)}\right)\left(\mathbf{r}\right)=\phi_{j}^{\left(m\right)}\left(\mathbf{r}\right)\left(-4\pi\mathcal{G}_{0}*\left(\sum_{i=1}^{N}\left|\phi_{i}^{\left(m\right)}\left(\mathbf{r}\right)\right|^{2}\right)\right)
\]
and
\[
\left(\mathcal{K}^{\left(m\right)}\phi_{j}\right)\left(\mathbf{r}\right)=\sum_{i=1}^{N}\phi_{i}^{\left(m\right)}\left(\mathbf{r}\right)\left(-4\pi\mathcal{G}_{0}*\left(\phi_{i}^{\left(m\right)}\phi_{j}^{\left(m\right)}\right)\right)
\]
(complex conjugation is removed noting that in our calculations orbitals
are real). Next, we compute entries of the matrix 
\begin{equation}
H_{ij}^{\left(m\right)}=-\frac{1}{2}\langle\Delta\phi_{i}^{\left(m\right)},\phi_{j}^{\left(m\right)}\rangle+\langle\mathcal{V}_{tot}^{\left(m\right)}\phi_{i}^{\left(m\right)},\phi_{j}^{\left(m\right)}\rangle,\,\,\,\,i,j=1,\dots N,\label{eq:HF iteration  matrix H}
\end{equation}
and evaluate its eigenvalues, which we set to be approximations of
the orbital energies $E_{i}^{\left(m\right)}$, $i=1,\dots N$. We
then compute 
\begin{equation}
\mu_{i}^{\left(m\right)}=\sqrt{-2E_{i}^{\left(m\right)}}\label{eq:HF iteration mu}
\end{equation}
and a new set of functions
\begin{eqnarray}
\begin{pmatrix}\tilde{\phi}_{1}^{\left(m+1\right)}\\
\tilde{\phi}_{2}^{\left(m+1\right)}\\
\vdots\\
\tilde{\phi}_{N}^{\left(m+1\right)}
\end{pmatrix} & \leftarrow-2 & \begin{pmatrix}\mathcal{G}_{\mu_{1}^{\left(m\right)}}*\left(\mathcal{V}_{tot}^{\left(m\right)}\phi_{1}^{\left(m\right)}\right)\\
\mathcal{G}_{\mu_{2}^{\left(m\right)}}*\left(\mathcal{V}_{tot}^{\left(m\right)}\phi_{2}^{\left(m\right)}\right)\\
\vdots\\
\mathcal{G}_{\mu_{N}^{\left(m\right)}}*\left(\mathcal{V}_{tot}^{\left(m\right)}\phi_{N}^{\left(m\right)}\right)
\end{pmatrix}.\label{eq:HF-iteration final orbital update}
\end{eqnarray}
Finally, we orthogonalize (and normalize) the resulting functions
$\tilde{\phi}_{j}^{\left(m+1\right)}$, $j=1,\dots,N$, yielding the
updated orbitals $\phi_{j}^{\left(m+1\right)}$, $j=1,\dots,N$.

We observe that as long as the Green's functions, the Coulomb potentials
and the orbitals are represented as Gaussian mixtures, then each operation
described above maintains this form, i.e. the result is again a Gaussian
mixture. However the number of terms grows rapidly and we show how
to control the number of terms in Section~\ref{sec:Reduction-algorithm}.
\begin{rem}
This iteration is known to be convergent for computing bound states
(see discussion in \cite{H-F-Y-G-B:2004}) although we are not aware
of a rigorous mathematical proof of this fact. An argument can be
made that if the potential $\mathcal{V}_{tot}$ is in the so-called
Rollnik class (see \cite{SIMON:1971}) then, for any fixed $-\mu_{i}^{2}<2E_{i}$,
$i=1,\dots,N$, the norms of relevant operators are less than $1$.
This makes this iteration a fixed point iteration. The nuclear Coulomb
potentials in $\mathcal{V}_{ext}$ just miss being in the Rollnik
class because of their slow decay away from the nuclei. However, any
truncation of these potentials that is sufficient for computing the
bound states (at an arbitrarily large distance from nuclei) puts them
back into the Rollnik class. In this paper we simply accept the fact
that the iteration is convergent. We do truncate the nuclear Coulomb
potentials at infinity in our approximations (see below) as is done
(explicitly or implicitly) in all numerical methods for computing
bound states.
\end{rem}

\subsection{\label{sec:Accurate-approximation-of-Green's}Accurate approximation
of Green's functions and potentials}

We briefly recall the key approximation results in \cite{BEY-MON:2005,BEY-MON:2010}
for the power functions $r^{-\alpha}$, $\alpha>0$, and the Green's
function for the non-oscillatory Helmholtz equations,
\[
\mathcal{G}_{\mu}\left(r\right)=\frac{1}{\left(2\pi\right)^{3/2}}\left(\frac{\mu}{r}\right)^{1/2}K_{\frac{1}{2}}\left(\mu r\right)=\frac{1}{4\pi}\frac{e^{-\mu r}}{r},
\]
where $K_{\frac{1}{2}}$ is the modified Bessel function of the second
kind. 
\begin{thm}
\cite[Theorem 5]{BEY-MON:2010} For any $\alpha>0$, $\delta>0$ and
$1/e\geq\epsilon>0$, there exist a step size $h$ and integers $M$
and $N$ such that
\[
\left|r^{-\alpha}-\mathcal{G}_{0}\left(r;M,N,h,\text{\ensuremath{\tau}}\right)\right|\leq r^{-\alpha}\epsilon,\ \ \ \mbox{for all}\ \delta\leq r\leq R,
\]
where
\begin{equation}
\mathcal{G}_{0}\left(r;M,N,h,\text{\ensuremath{\tau}}\right)=\frac{h}{\Gamma\left(\alpha/2\right)}\sum_{n=M+1}^{N}e^{\alpha\left(hn-\tau\right)/2}e^{-e^{hn-\tau}r^{2}},\label{eq:est in bey-mon-2010}
\end{equation}
and $\tau$ is any number \textup{$0\le\tau<h$.}
\end{thm}

The error estimates are based on discretizing the integral
\[
\frac{1}{r^{\alpha}}=\frac{1}{\Gamma\left(\frac{\alpha}{2}\right)}\int_{-\infty}^{\infty}e^{-r^{2}e^{t}+\frac{\alpha}{2}t}dt
\]
using (an infinite) trapezoidal rule
\begin{equation}
\frac{1}{\Gamma\left(\frac{\alpha}{2}\right)}\int_{-\infty}^{\infty}e^{-r^{2}e^{t}+\frac{\alpha}{2}t}dt\approx\frac{h}{\Gamma\left(\frac{\alpha}{2}\right)}\sum_{n\in\mathbb{Z}}e^{\alpha\left(hn-\tau\right)/2}e^{-e^{hn-\tau}r^{2}},\label{eq:infinite sum}
\end{equation}
where the step size $h$ satisfies
\[
h\leq\frac{2\pi}{\log3+\frac{\alpha}{2}\log(\cos1)^{-1}+\log\epsilon^{-1}}
\]
and $\epsilon$ is any user-selected accuracy. The proof in \cite[Theorem 5]{BEY-MON:2010}
allows one to shift the grid by any $\tau$, $0\le\tau<h$.

For a given accuracy $\epsilon$, power $\alpha$ and a range of values
$r$, the infinite sum (\ref{eq:infinite sum}) is then truncated
due to the exponential or super-exponential decay of the integrand
at $\pm\infty$, to yield a finite sum approximation in that range.
The number of terms is estimated as
\[
N-M\leq\frac{1}{10}(2\log\epsilon^{-1}+\log\alpha+2)(\log\delta^{-1}+\frac{1}{\alpha}\log\epsilon^{-1}+\log\log\epsilon^{-1}+\frac{3}{2}).
\]
This approximation provides an analytic construction of a multiresolution,
separated representation for the Poisson kernel and the Coulomb potential
in any dimensions. In practice many terms with small exponents in
(\ref{eq:est in bey-mon-2010}) can be combined using \cite[Section 6]{BEY-MON:2005}
to reduce the total number of terms further. In our examples we use
the following approximation of the Poisson kernel and the Coulomb
potential via a linear combination of $146$ Gaussians:
\begin{equation}
\left|r^{-1}-\sum_{n=1}^{8}w_{n}e^{-\eta_{n}r^{2}}-\frac{h}{\Gamma\left(\frac{1}{2}\right)}\sum_{n=-51}^{87}e^{\left(hn-\tau\right)/2}e^{-e^{hn-\tau}r^{2}}\right|\leq r^{-1}\epsilon,\ \ \ \mbox{for all}\ 10^{-7}\leq r\leq10^{5},\label{eq:approx of 1/r}
\end{equation}
where $h=0.40994422603935795$, $\tau=0.192967891816239$, $\epsilon=10^{-10}$
and exponents $\eta_{n}$ and coefficients $w_{n}$ are given in Table~\ref{tab:Parameters-in-}.
\begin{center}
\begin{table}
\begin{centering}
\begin{tabular}{|c|c|c|}
\hline 
$n$ & $\eta_{n}$ & $w_{n}$\tabularnewline
\hline 
\hline 
1 & $2.1073876854180\cdot10^{-12}$ & $3.2630674210379\cdot10^{-6}$\tabularnewline
\hline 
2 & $1.8365780986634\cdot10^{-11}$ & $3.1058837221013\cdot10^{-6}$\tabularnewline
\hline 
3 & $4.7777245228151\cdot10^{-11}$ & $2.8014247111005\cdot10^{-6}$\tabularnewline
\hline 
4 & $8.5624630300630\cdot10^{-11}$ & $2.5227064974618\cdot10^{-6}$\tabularnewline
\hline 
5 & $1.3289239111902\cdot10^{-10}$ & $2.7039982943831\cdot10^{-6}$\tabularnewline
\hline 
6 & $2.0054640049463\cdot10^{-10}$ & $3.2761422288967\cdot10^{-6}$\tabularnewline
\hline 
7 & $3.0217586807074\cdot10^{-10}$ & $4.0205002817225\cdot10^{-6}$\tabularnewline
\hline 
8 & $4.5529860118663\cdot10^{-10}$ & $4.9351231646262\cdot10^{-6}$\tabularnewline
\hline 
\multicolumn{1}{c}{} & \multicolumn{1}{c}{} & \multicolumn{1}{c}{}\tabularnewline
\end{tabular}
\par\end{centering}
\caption{Exponents and weights for the first eight terms in (\ref{eq:approx of 1/r}).
\label{tab:Parameters-in-} }
\end{table}
\par\end{center}

A similar result holds for the Green's function for the non-oscillatory
Helmholtz equation,
\begin{equation}
\mathcal{G}_{\mu}\left(r\right)=\frac{1}{4\pi}\frac{e^{-\mu r}}{r}=\left(4\pi\right)^{-\frac{3}{2}}\int_{-\infty}^{\infty}e^{-\frac{r^{2}e^{t}}{4}-\mu^{2}e^{-t}+\frac{t}{2}}dt,\label{eq:GF with mu}
\end{equation}
where the integrand has a super-exponential decay at $\pm\infty$.
An accurate and efficient approximation is obtained by discretizing
this integral representation. We use 
\[
\mathcal{G}_{\mu}\left(r;M,N,h\right)=\left(4\pi\right)^{-3/2}h\sum_{l=M}^{N}e^{-\mu^{2}e^{-hl}+\frac{hl}{2}}e^{-\frac{\left\Vert \mathbf{r}\right\Vert ^{2}e^{hl}}{4}}
\]
where the step size $h$ is selected to achieve the desired accuracy
$\epsilon$, 
\[
\left|\mathcal{G}_{\mu}\left(r\right)-\mathcal{G}_{\mu}\left(r;M,N,h\right)\right|\leq r^{-1}\epsilon,\ \ \ \mbox{for all}\ \delta\leq r\leq R.
\]
Note that this integral is discretized for a range of the parameter
$\mu$ so that only the weights (coefficients) change as we change
$\mu$. In our examples we approximate the Green's function $\mathcal{G}_{\mu}\left(r\right)$
(\ref{eq:GF with mu}) with 
\begin{equation}
\mathcal{G}_{\mu}\left(r;-20,120,0.38190954773869346734\right)\label{eq:approx G_mu}
\end{equation}
resulting in $\epsilon=10^{-10}$ for $r$ in the range $10^{-7}\leq r\leq10^{5}$. 

\section{\label{sec:Reduction-algorithm}Reduction algorithm with subdivision
of terms into groups}

We approximate orbitals using Gaussian mixtures (cf. the multiwavelets
in \cite{H-F-Y-G-B:2004,Y-F-G-H-B:2004,Y-F-G-H-B:2004a} and the Slater-type
orbitals in \cite{BEY-HAU:2013}). While using a fixed basis set of
Gaussians selected in advance has a long history in quantum chemistry,
in contrast with previous methods we do not choose a set of Gaussians
in advance (nor do we use a particular multiresolution basis of Gaussians
although such an approach is possible using results in \cite{BE-MO-SA:2017}).
We allow the iterative algorithm for solving the equations in Section~\ref{subsec:Iteration-for-solving}
to ``select'' the necessary basis functions. After each operation
that yields a significant increase in the number of Gaussians, we
prune the resulting Gaussian mixture by using the reduction algorithm
in \cite[Algorithm 1]{BE-MO-YA:2019} that chooses the ``best''
linearly independent subset of terms, the so-called skeleton terms.
By estimating the error of the solution, we terminate the iteration
when the desired accuracy is achieved. 

Conceptually our method is simpler and, also, less technical than
either \cite{H-F-Y-G-B:2004,Y-F-G-H-B:2004,Y-F-G-H-B:2004a} or \cite{BEY-HAU:2013}
and, as far as we know, is novel. Clearly, the key to enabling our
approach is the reduction algorithm described and analyzed in \cite{BE-MO-YA:2019}.
This reduction algorithm has complexity $\mathcal{O}\left(r^{2}N+p\left(d\right)rN\right)$,
where $N$ is the initial number of terms, $r$ is the number of skeleton
terms and $p\left(d\right)$ is the cost of computing the inner product
between the terms of a Gaussian mixture ($p\left(d\right)$ is a small
constant in dimension $d=3$). If the number of skeleton terms, $r$,
is small, the computational cost of the reduction algorithm is completely
satisfactory. However, we cannot assume that the number of skeleton
terms (which are used to represent orbitals) is small for a large
molecule. 

In this paper we introduce an important modification of the reduction
algorithm to avoid an excessive computational cost due to the quadratic
dependence on the number of skeleton terms. The modification stems
from a simple observation that, if two sets of functions are known
to be linearly independent in advance, an attempt to find skeleton
terms is wasteful. In many applied problems a multivariate mixture
with a large number of terms can be split into subgroups of terms
so that each subgroup reflects the local behavior of the function.
For example, it is well understood that orbitals have cusps at nuclear
centers but otherwise are smooth functions (in fact it is this structure
of solutions that underpins the success of multiresolution methods
such as \cite{H-F-Y-G-B:2004,Y-F-G-H-B:2004,Y-F-G-H-B:2004a}). Well-localized
terms associated with cusps decay rapidly away from a nuclear center
so that these terms will be linearly independent (or even nearly orthogonal)
from similar terms at other nuclear centers. This suggests an acceleration
technique for the reduction algorithm based on a hierarchical subdivision
of the terms of the mixture by scale and location. Given such a subdivision,
reduction is performed only within a local group. In such an approach
the reduction of terms within each group is independent from other
groups and, thus, the computation is trivially parallel. 

In \cite{BE-MO-YA:2019} we consider Gaussian terms in the form of
$L^{2}$-normalized multivariate normal distribution $N\left(\mathbf{x},\mathbf{s}_{k},\boldsymbol{\Sigma}_{k}\right)$,
\begin{eqnarray}
g_{k}\left(\mathbf{x},\mathbf{s}_{k},\boldsymbol{\Sigma}_{k}\right) & = & \left(\det\left(4\pi\boldsymbol{\Sigma}_{k}\right)\right)^{\frac{1}{4}}N\left(\mathbf{x},\mathbf{s}_{k},\boldsymbol{\Sigma}_{k}\right)\label{eq:Gaussian_atoms}\\
 & = & \frac{1}{\left(\det\left(\pi\boldsymbol{\Sigma}_{k}\right)\right)^{1/4}}\exp\left(-\frac{1}{2}\left(\mathbf{x}-\mathbf{s}_{k}\right)^{T}\boldsymbol{\Sigma}_{k}^{-1}\left(\mathbf{x}-\mathbf{s}_{k}\right)\right),\nonumber 
\end{eqnarray}
where the vector $\mathbf{s}_{k}$ defines its location and the symmetric
positive definite matrix $\boldsymbol{\Sigma}_{k}$ controls it shape.
In the setup of this paper, all matrices $\boldsymbol{\Sigma}_{k}$
are diagonal and, therefore, functions represented as Gaussian mixtures
are in separated form (see \cite{BEY-MOH:2002,BEY-MOH:2005}). In
fact, in our examples it is sufficient to use rotationally symmetric
Gaussians so that their shape is controlled by a single scalar parameter.

We subdivide terms by grouping them by their scale (shape) and location.
We consider a single (global) group of flat Gaussians. i.e. Gaussians
that have a significant support around more than one nuclear center
and control behavior of orbitals far away from nuclear centers. The
rest of the terms we split into groups according to their proximity
to nuclear centers using nuclear centers as ``seeds'' in a Voronoi-type
decomposition: all terms are split into groups by their proximity
to the nuclei, i.e.\ a term located at $\mathbf{s}_{k}$ belongs
to a group associated with a given nucleus at location $\mathbf{R}_{l}$
if it is the closest to it among all other nuclei, i.e. $\left\Vert \mathbf{R}_{l}-\mathbf{s}_{k}\right\Vert \le\min_{l'}\left\Vert \mathbf{R}_{l'}-\mathbf{s}_{k}\right\Vert $.
Such subdivision ensures that the cost of reduction is proportional
to the number of nuclei since the number of skeleton terms associated
with each group is expected to be roughly the same given the structure
of functions they represent. In our experiments we observed that the
number of skeleton terms in these groups is still fairly large (e.g.
several thousands terms) so further subdivision by scale and location
is appropriate. Towards this end, we associate with each nucleus center
and its group the distance $s_{l}^{\max}=\max_{k}\left\Vert \mathbf{R}_{l}-\mathbf{s}_{k}\right\Vert $
and subdivide each group further by considering two parameters: (i)
the radial intervals $\left[\frac{m}{2^{j}}s_{l}^{\max},\frac{m+1}{2^{j}}s_{l}^{\max}\right]$,
$m=0,1,\dots,2^{j}-1$ that restrict the locations of the terms; (ii)
the ranges of exponents $\left[4^{-j-1}\sigma^{\mbox{far}},4^{-j}\sigma^{\mbox{far}}\right]$,
where $\sigma^{\mbox{far}}$ is a parameter (see discussion below),
identifying the shape of the terms. Within each radial interval $\left[\frac{m}{2^{j}}s_{l}^{\max},\frac{m+1}{2^{j}}s_{l}^{\max}\right]$,
terms may have significantly different shapes determined by the shape
matrices $\Sigma_{k}$ and, thus, cannot be linearly dependent. Effectively,
it is a multiresolution subdivision strategy. 

In the examples in Section~\ref{sec:Examples}, we initialize the
spatial orbitals so that every term in the mixture has a diagonal
shape matrix $\Sigma_{k}=\mbox{diag}\left(\sigma_{k},\dots,\sigma_{k}\right)=\sigma_{k}I$,
where $I$ is the identity matrix. This forces all Gaussians involved
in further computations to be rotationally symmetric, i.e. to have
$\Sigma_{k}$ of the same structure. As already stated, we first identify
a global group by considering all terms with $\sigma_{k}\in\left[\sigma^{\mbox{far}},\infty\right)$,
where $\sigma^{\mbox{far}}>0$ is a parameter. Terms assigned in this
group are significant in the vicinity of more than one nuclear center
and capture the far-field behavior. The rest of the terms we divide
into groups by their proximity to the nuclei and scale (shape). Within
each group associated with nuclei $\mathbf{R}_{l}$, we identify a
Gaussian term $g_{k}$ as a member of a subgroup, denoted by $G_{jm}^{l}$
if 
\[
\left\Vert \mathbf{s}_{k}-\mathbf{s}_{l}^{\mbox{max}}\right\Vert \in\left[\frac{m}{2^{j}}s_{l}^{\max},\frac{m+1}{2^{j}}s_{l}^{\max}\right]\ \ \mbox{and}\ \ \sigma_{k}\in\left[4^{-j-1}\sigma^{\mbox{far}},4^{-j}\sigma^{\mbox{far}}\right],
\]
where $j=0,1,\dots,J$. The index $m$ controls the location and the
index $j$ controls the scale (shape) of a term. The role of these
parameters is similar to those labeling basis functions in Multiresolution
Analysis (MRA).
\begin{rem}
In general, subdivision into groups may differ in detail from the
one we use for our examples in this paper. This depends on the type
of functions we want to represent. In particular, the number of indices
controlling location of terms (and, therefore, subgroups) can be larger
than what we use in our examples (e.g. there could be three such indices).
\end{rem}

We note that, for large $j$ and $m$, the corresponding group $G_{jm}^{l}$
consists of terms with a small essential support and located (relatively)
far away from the nucleus center $\mathbf{R}_{l}$. Such terms (and,
therefore, groups) are not needed to represent orbitals since orbitals
have cusps only at the nuclei. To reduce the number of unnecessary
groups, we do not maintain such groups. 

In order to capture the cusp behavior, we introduce in the vicinity
of nuclei an additional set of groups $G_{j,0}^{l}$, $j=J+1,\dots,\widetilde{J},$
so that $g_{k}\in G_{j,0}^{l}$ if
\[
\left\Vert \mathbf{s}_{k}-\mathbf{s}_{l}^{\mbox{max}}\right\Vert \in\left[0,\frac{s_{l}^{\mbox{max}}}{2^{j}}\right]\ \ \mbox{and}\ \ \sigma_{k}\in\left[4^{-j-1}\sigma^{\mbox{far}},4^{-j}\sigma^{\mbox{far}}\right].
\]
In our computation we ignore any terms that are not included in the
groups $\left\{ G_{jm}^{l}\right\} _{j=0,\dots,J}^{m=0,\dots,2^{j}-1}$
or $\left\{ G_{j,0}^{l}\right\} _{j=J+1,\dots,\widetilde{J}}$ for
all nucleus locations $\mathbf{R}_{l}$. Selecting $J=4$ and $\widetilde{J}=26$,
the total number of groups (including the global group) is $107$.
It is worth noting that, when reducing a Gaussian mixture according
to the group subdivision, some of the groups might be empty.

To estimate how subdivision into groups speeds up computations, let
us consider a multivariate mixture with $N$ terms where the number
of linearly independent terms is $r$. Let us assume that we can subdivide
these terms into $m$ groups with an equal number of terms (for simplicity
of the estimate). We assume that the number of linearly independent
terms in each group is then at most $r/k$, where $k>1$. The cost
of performing the reduction algorithm for one group is then $\mathcal{O}\left(\left(\frac{r}{k}\right)^{2}\frac{N}{m}+p\left(d\right)\frac{r}{k}\frac{N}{m}\right)$
and, thus, the total cost to reduce all groups is $\mathcal{O}\left(\left(\frac{r}{k}\right)^{2}N+p\left(d\right)\frac{r}{k}N\right)$.
Depending on $k$, we get a significant speed-up factor and note that
the reduction of independent groups of terms is trivially parallel
(a property that we did not use it in our computations). With the
subdivision into groups described above, the number of groups is naturally
proportional to the number of nuclear centers. Since the solution
has a cusp at the nuclear center, we can consider $r/k$ to be a constant
(i.e., each cusp locally requires roughly the same number of terms).
Hence, the reduction algorithm with subdivision into groups is linear
in the number of nuclear centers. 
\begin{rem}
\label{rem:By-examining-operators}By examining operators in (\ref{eq:HF-iteration})-(\ref{eq:HF-iteration final orbital update}),
we observe that if we start an iteration with a Gaussian mixture with
diagonal shape matrices $\Sigma_{k}$, then all resulting Gaussian
mixtures involved in the computation will have diagonal shape matrices
as well. Moreover, location of the terms is restricted to a box defined
by positions of nuclei. We observe that, in (\ref{eq:HF-iteration})-(\ref{eq:HF-iteration final orbital update}),
convolutions and multiplications produce new Gaussians (with different
shifts and/or shape matrices). While convolution does not change the
center of the resulting Gaussian terms, multiplication of two Gaussians
generates a Gaussian centered at a different location. Specifically,
we have for the product of two normal distributions

\begin{equation}
N\left(\mathbf{x},\boldsymbol{\mu}_{1},\boldsymbol{\Sigma}_{1}\right)N\left(\mathbf{x},\boldsymbol{\mu}_{2},\boldsymbol{\Sigma}_{2}\right)=N\left(\boldsymbol{\mu}_{1},\boldsymbol{\mu}_{2},\boldsymbol{\Sigma}_{1}+\boldsymbol{\Sigma}_{2}\right)\cdot N\left(\mathbf{x},\boldsymbol{\mu}_{c},\left(\boldsymbol{\Sigma}_{1}^{-1}+\boldsymbol{\Sigma}_{2}^{-1}\right)^{-1}\right)\label{eq:product}
\end{equation}
where
\[
\boldsymbol{\mu}_{c}=\left(\boldsymbol{\Sigma}_{1}^{-1}+\boldsymbol{\Sigma}_{2}^{-1}\right)^{-1}\left(\boldsymbol{\Sigma}_{1}^{-1}\boldsymbol{\mu}_{1}+\boldsymbol{\Sigma}_{2}^{-1}\boldsymbol{\mu}_{2}\right).
\]
If $\boldsymbol{\Sigma}_{1}=\mbox{diag}\left(\alpha_{1},\alpha_{2},\alpha_{3}\right)$,
$\boldsymbol{\Sigma}_{2}=\mbox{diag}\left(\beta_{1},\beta_{2},\beta_{3}\right)$,
$\boldsymbol{\mu}_{1}=\left(x_{1},x_{2},x_{3}\right)$ and $\boldsymbol{\mu}_{2}=\left(y_{1},y_{2},y_{3}\right)$,
then we have 
\[
\boldsymbol{\mu}_{c}=\left(\frac{\beta_{1}}{\alpha_{1}+\beta_{1}}x_{1}+\frac{\alpha_{1}}{\alpha_{1}+\beta_{1}}y_{1},\frac{\beta_{2}}{\alpha_{2}+\beta_{2}}x_{2}+\frac{\alpha_{2}}{\alpha_{2}+\beta_{2}}y_{2},\frac{\beta_{3}}{\alpha_{3}+\beta_{3}}x_{3}+\frac{\alpha_{3}}{\alpha_{3}+\beta_{3}}y_{3}\right).
\]
Hence, all three coordinates of the new center $\boldsymbol{\mu}_{c}$
are in between the corresponding coordinates of $\boldsymbol{\mu}_{1}$
and $\boldsymbol{\mu}_{2}$. This implies that all Gaussian terms
resulting from iteration (\ref{eq:HF-iteration})-(\ref{eq:HF-iteration final orbital update})
will have their centers inside a box defined by the minimum and the
maximum of individual coordinates of nuclear centers $\mathbf{R}_{l}$
in (\ref{eq: external potential}) (in fact, this result holds in
any dimension). Moreover, if the Gaussians are rotationally invariant,
then all terms will have their centers inside a convex polyhedron
defined by the nuclear centers $\mathbf{R}_{l}$. We always initialize
the iteration using Gaussians centered inside the convex polyhedron
defined by the nuclear centers.

We note that the restricted location of terms of a Gaussian mixture
in our approach cannot be achieved when using wavelet or multiwavelet
bases, which force a computational box that is much larger than the
box defined by nuclear centers.
\end{rem}

\section{\label{sec:Examples}Examples of computations for two small molecules}

\subsection{\label{subsec:hellium}Helium Hydride, $\mbox{HeH}^{+}$}

First we solve the Hartree\textendash Fock equation for $\mbox{HeH}^{+}$
\begin{equation}
\left(-\frac{1}{2}\Delta+\mathcal{V}_{ext}-4\pi\Delta^{-1}\left(\left|\phi\right|^{2}\right)\right)\phi=E\phi,\label{eq:helium eqn}
\end{equation}
with the potential
\[
\mathcal{V}_{ext}\left(\mathbf{r}\right)=\frac{Z_{1}}{\left\Vert \mathbf{r}-\mathbf{R}_{1}\right\Vert }+\frac{Z_{2}}{\left\Vert \mathbf{r}-\mathbf{R}_{2}\right\Vert }.
\]
where $Z_{1}=-1$, $Z_{2}=-2$, $\mathbf{R}_{1}=\left(0,0,-0.7\right)$
and $\mathbf{R}_{2}=\left(0,0,0.7\right)$. 

As described in Section~\ref{subsec:Iteration-for-solving}, our
basic approach involves recasting (\ref{eq:helium eqn}) as an integral
equation which we solve iteratively. The iteration (see below) is
convergent and, while the rate of convergence is slower than quadratic,
it is sufficiently fast so that only a few dozen iterations are required
(see discussion in \cite{H-F-Y-G-B:2004} and references therein).
We represent the spatial orbital $\phi\left(\mathbf{r}\right)$ as
a Gaussian mixture
\[
\phi\left(\mathbf{r}\right)=\sum_{k=1}^{K}c_{k}g_{k}\left(\mathbf{r},\mathbf{s}_{k},\boldsymbol{\Sigma}_{k}\right),\ \ \mathbf{r}\in\mathbb{R}^{3},
\]
that is gradually constructed via iterations. In each step of the
iterative solution, the number of terms in representing $\phi\left(\mathbf{r}\right)$
grows and we use \cite[Algorithm 1]{BE-MO-YA:2019} to control their
number.

We approximate potentials $Z_{1}/\left\Vert \mathbf{r}-\mathbf{R}_{1}\right\Vert $
and $Z_{2}/\left\Vert \mathbf{r}-\mathbf{R}_{2}\right\Vert $ and
the Poisson kernel using (\ref{eq:approx of 1/r}) and the Green's
function $\mathcal{G}_{\mu}\left(r\right)$ in (\ref{eq:GF with mu})
using (\ref{eq:approx G_mu}). 

For the reader's convenience, we explicitly describe the iteration
for solving (\ref{eq:helium eqn}) ,
\begin{eqnarray}
\mathcal{V}_{tot}^{\left(m\right)}\phi^{\left(m\right)} & \leftarrow & \left(\mathcal{V}_{ext}-4\pi\mathcal{G}_{0}*\left(\left|\phi^{\left(m\right)}\right|^{2}\right)\right)\phi^{\left(m\right)}\label{eq:helium iteration}\\
E^{\left(m\right)} & \leftarrow & \left\langle -\frac{1}{2}\Delta\phi^{\left(m\right)}+\mathcal{V}_{tot}^{\left(m\right)}\phi^{\left(m\right)},\phi^{\left(m\right)}\right\rangle ,\nonumber \\
\mu^{\left(m\right)} & \leftarrow & \sqrt{-2E^{\left(m\right)}},\nonumber \\
\phi^{\left(m+1\right)} & \leftarrow & -2\mathcal{G}_{\mu^{\left(m\right)}}*\left(\mathcal{V}_{tot}^{\left(m\right)}\phi^{\left(m\right)}\right),\nonumber \\
\phi^{\left(m+1\right)} & \leftarrow & \frac{\phi^{\left(m+1\right)}}{\left\Vert \phi^{\left(m+1\right)}\right\Vert _{2}},\ \ \ m=0,1,\cdots,\nonumber 
\end{eqnarray}
where $*$ denotes convolution. 

To start the iteration, we use the single Gaussian 
\[
\phi^{\left(0\right)}\left(\mathbf{r}\right)=\pi^{-\frac{3}{4}}e^{-\frac{\left\Vert \mathbf{r}\right\Vert ^{2}}{2}}.
\]
We choose $\epsilon=10^{-6}$ as the error tolerance for the reduction
of terms in each group and terminate the iteration if $\left|E^{\left(m\right)}-E^{\left(m-1\right)}\right|<4\times10^{-6}$.
After $20$ iterations we obtain the orbital energy $E=-1.6605545$
and the total energy $E_{tot}=-2.9325704$ . Using the MADNESS software
\cite{F-H-B-H-S-S:2007} yields orbital energy $E=-1.66053903$ and
total energy $-2.93256741$ so that, using these values as reference,
our computations have absolute error $1.5\cdot10^{-5}$ for the orbital
energy and $3.0\cdot10^{-6}$ for the total energy. As a result of
the iteration, the total number of terms representing the orbital
$\phi\left(\mathbf{r}\right)$ is $1563$. This number of terms is
significantly smaller than the $\approx2\cdot10^{5}$ terms required
when using multiwavelets \cite{H-F-Y-G-B:2004,Y-F-G-H-B:2004,Y-F-G-H-B:2004a}
and a factor $\approx2$ larger than that in \cite{BEY-HAU:2013},
where the number of Slater-type orbitals is $637$. However, our approach
is less technical than that in \cite{BEY-HAU:2013} and, for this
reason, it is easier to extend it to deal with large molecules. In
Table~\ref{tab:HeH}, we show the number of terms in the global group
$N_{global}$, the number of non-empty groups $N_{groups}$ and the
maximum $N_{max}$, minimum $N_{min}$ and average $N_{ave}$ number
of terms in these non-empty groups in $\phi$ and $\mathcal{V}_{tot}\phi$
at the final step.
\begin{table}
\centering{}%
\begin{tabular}{|c|c|c|c|c|c|}
\hline 
 & $N_{global}$ & $N_{groups}$ & $N_{max}$ & $N_{min}$ & $N_{ave}$\tabularnewline
\hline 
\hline 
$\phi$ & $137$ & $28$ & $69$ & $2$ & $50.9$\tabularnewline
\hline 
$\mathcal{V}_{tot}\phi$ & $137$ & $70$ & $70$ & $4$ & $32.8$\tabularnewline
\hline 
\multicolumn{1}{c}{} & \multicolumn{1}{c}{} & \multicolumn{1}{c}{} & \multicolumn{1}{c}{} & \multicolumn{1}{c}{} & \multicolumn{1}{c}{}\tabularnewline
\end{tabular}\caption{\label{tab:HeH}The number of terms in the global group $N_{global}$,
the number of non-empty groups $N_{groups}$ and the maximum $N_{max}$,
minimum $N_{min}$ and average $N_{ave}$ number of terms in these
non-empty groups in Gaussian mixtures representing $\phi$ and $\mathcal{V}_{tot}\phi$
at the final step.}
\end{table}

\begin{figure}[h]
\begin{centering}
\includegraphics[scale=0.5]{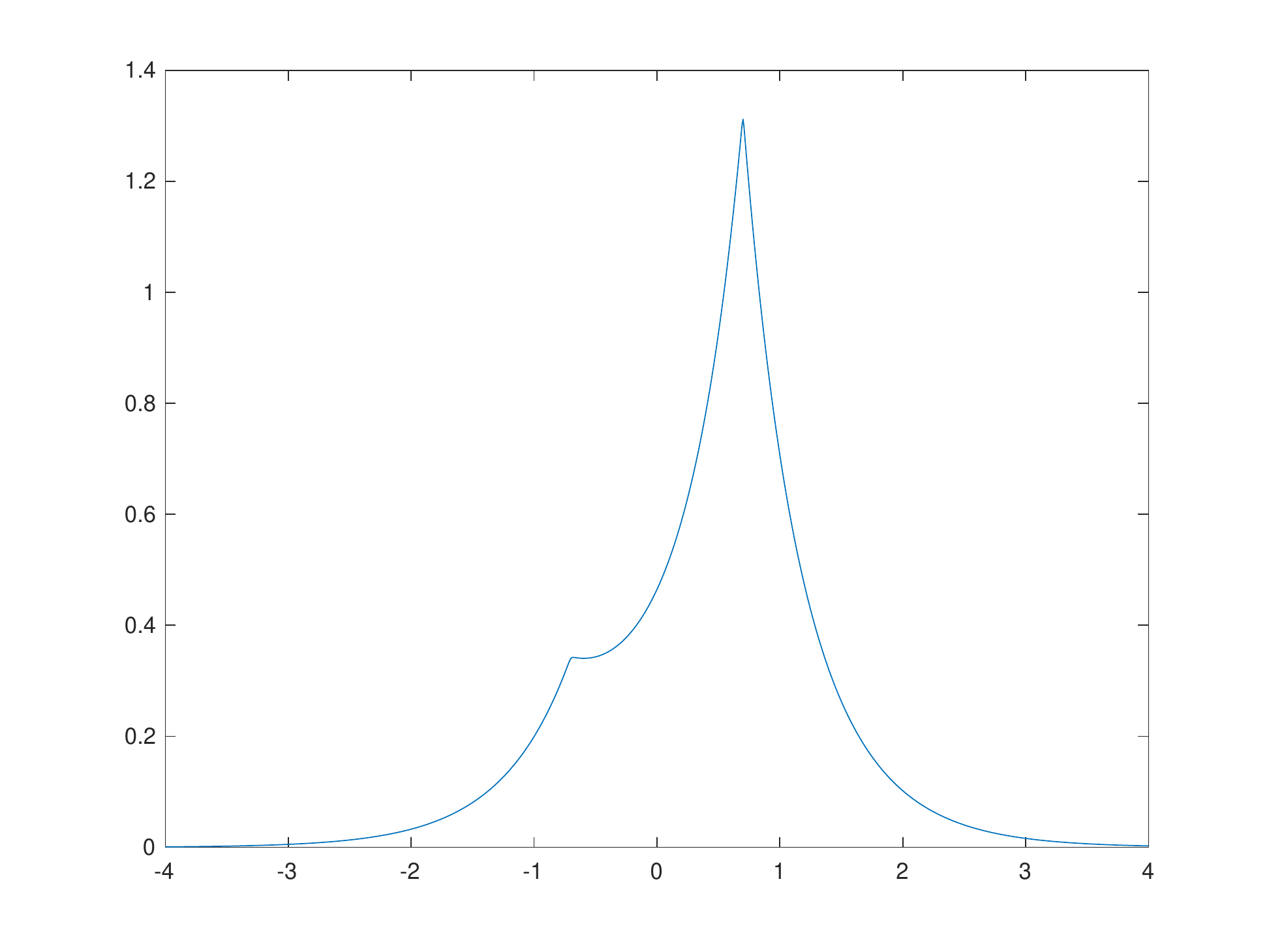}
\par\end{centering}
\caption{\label{fig:helium}Plot of the spatial orbital $\phi\left(\mathbf{r}\right)$
for helium hydride $\mbox{HeH}^{+}$ on the line $\mathbf{r}=\left(0,0,x\right)$
connecting the two nuclear centers $\mathbf{R}_{1}=\left(0,0,-0.7\right)$
and $\mathbf{R}_{2}=\left(0,0,0.7\right)$.}
\end{figure}

In Figure\,\ref{fig:helium}, we display the spatial orbital $\phi\left(\mathbf{r}\right)$
on the line $\mathbf{r}=\left(0,0,x\right)$ connecting the two nuclei
locations $\mathbf{R}_{1}$ and $\mathbf{R}_{2}$. 

\subsection{\label{subsec:LiH}Lithium Hydride, LiH}

For our second example, we consider the Hartree-Fock equations for
lithium hydride, $\mbox{LiH}$. The Hartree\textendash Fock equations
in this case are 
\begin{equation}
\mathcal{F}\phi_{j}\left(\mathbf{r}\right)=E_{j}\phi\left(\mathbf{r}\right),\ \ j=1,2,\label{eq:lithium eqn}
\end{equation}
where $\mathcal{F}=-\frac{1}{2}\Delta+\mathcal{V}_{ext}+2\mathcal{J}-\mathcal{K},$
\[
\mathcal{J}\phi_{j}=\phi_{j}\left(-4\pi\Delta^{-1}\left(\left|\phi_{1}\right|^{2}+\left|\phi_{2}\right|^{2}\right)\right),
\]
\[
\mathcal{K}\phi_{j}=\phi_{1}\left(-4\pi\Delta^{-1}\left(\phi_{1}^{*}\phi_{j}\right)\right)+\phi_{2}\left(-4\pi\Delta^{-1}\left(\phi_{2}^{*}\phi_{j}\right)\right),
\]
 and 
\[
\mathcal{V}_{ext}\left(\mathbf{r}\right)=\frac{Z_{1}}{\left\Vert \mathbf{r}-\mathbf{R}_{1}\right\Vert }+\frac{Z_{2}}{\left\Vert \mathbf{r}-\mathbf{R}_{2}\right\Vert }.
\]
For lithium hydride we have $Z_{1}=-3$, $Z_{2}=-1$ , $\mathbf{R}_{1}=\left(-3.15/2,0,0\right)$
and $\mathbf{R}_{2}=\left(3.15/2,0,0\right)$. We approximate potentials
$Z_{1}/\left\Vert \mathbf{r}-\mathbf{R}_{1}\right\Vert $ and $Z_{2}/\left\Vert \mathbf{r}-\mathbf{R}_{2}\right\Vert $
and the Poisson kernel using (\ref{eq:approx of 1/r}) and the Green's
function $\mathcal{G}_{\mu}\left(r\right)$ in (\ref{eq:GF with mu})
using (\ref{eq:approx G_mu}). 

To solve (\ref{eq:lithium eqn}), we initialize orbitals $\phi_{j}^{\left(0\right)}\left(\mathbf{r}\right),\ j=1,2$
as
\[
\phi_{1}^{\left(0\right)}\left(\mathbf{r}\right)=e^{-\frac{\left\Vert \mathbf{r}-\mathbf{R}_{1}\right\Vert ^{2}}{20}}\ \ \mbox{and}\ \ \phi_{2}^{\left(0\right)}\left(\mathbf{r}\right)=e^{-\frac{\left\Vert \mathbf{r}-\mathbf{R}_{2}\right\Vert ^{2}}{20}}.
\]
There are many possible initializations, e.g., we can use an approximation
to $e^{-\left\Vert \mathbf{r}-\mathbf{R}_{1}\right\Vert }$ and $e^{-\left\Vert \mathbf{r}-\mathbf{R}_{2}\right\Vert }$
via Gaussian mixtures. Importantly, the initial approximation should
be chosen so that the initial orbital energies are negative. After
orthonormalizing the functions $\phi_{j}^{\left(0\right)}\left(\mathbf{r}\right)$,
we compute the initial energies $E_{j}^{\left(0\right)}$ as the eigenvalues
of the $2\times2$ matrix defined in (\ref{eq:HF iteration  matrix H}).
We then set $\mu_{j}^{\left(0\right)},\ j=1,2$ using (\ref{eq:HF iteration mu})
and update the orbitals via (\ref{eq:HF-iteration final orbital update}). 

Finally, when the desired accuracy is achieved, we compute the total
energy $E_{tot}$ as
\[
E_{tot}=\sum_{j=1}^{2}\left(E_{j}+\left\langle \left(-\frac{1}{2}\Delta+\mathcal{V}_{ext}\right)\phi_{j},\phi_{j}\right\rangle \right)+\frac{Z_{1}Z_{2}}{\left\Vert \boldsymbol{R}_{1}-\boldsymbol{R}_{2}\right\Vert }.
\]

In this example, we choose an error tolerance of $\epsilon=10^{-6}$
for the reduction of terms in each group. We terminate the computation
if $\left|E_{j}^{\left(m\right)}-E_{j}^{\left(m-1\right)}\right|<4\times10^{-6}$
for $j=1,2$. After $28$ iterations, we arrived at the orbital energies
$E_{1}=-2.451757267$ and $E_{2}=-0.297819313$ computed with absolute
errors of $5.7\times10^{-6}$ and $3.7\times10^{-6}$, respectively,
when compared with the reference energies $E_{1}=-2.451763$ and $E_{2}=-0.297823$
obtained using the MADNESS software. The computed total energy $E_{tot}=-7.9869324$
has an absolute error of $4.0\times10^{-6}$ compared with the value
$E_{tot}=-7.9869364$ evaluated in MADNESS. The total number of terms
to represent orbitals $\phi_{1}\left(\mathbf{r}\right)$ and $\phi_{2}\left(\mathbf{r}\right)$
is $2185$ and $2569$. In Table~\ref{tab:LiH}, we show the number
of terms in the global group $N_{global}$, the number of non-empty
groups $N_{groups}$ and the maximum $N_{max}$, minimum $N_{min}$
and average $N_{ave}$ number of terms in the non-empty groups in
the Gaussian mixture representations of $\phi_{1}$, $\phi_{2}$,
$\mathcal{V}_{tot}\phi_{1}$ and $\mathcal{V}_{tot}\phi_{2}$ at the
final step.
\begin{table}
\begin{centering}
\begin{tabular}{|c|c|c|c|c|c|}
\hline 
 & $N_{global}$ & $N_{group}$ & $N_{max}$ & $N_{min}$ & $N_{ave}$\tabularnewline
\hline 
\hline 
$\phi_{1}$ & $218$ & $26$ & $108$ & $1$ & $75.7$\tabularnewline
\hline 
$\phi_{2}$ & $217$ & $31$ & $107$ & $1$ & $75.9$\tabularnewline
\hline 
$\mathcal{V}_{tot}\phi_{1}$ & $220$ & $68$ & $107$ & $3$ & $43.0$\tabularnewline
\hline 
$\mathcal{V}_{tot}\phi_{2}$ & $215$ & $69$ & $107$ & $4$ & $46.2$\tabularnewline
\hline 
\multicolumn{1}{c}{} & \multicolumn{1}{c}{} & \multicolumn{1}{c}{} & \multicolumn{1}{c}{} & \multicolumn{1}{c}{} & \multicolumn{1}{c}{}\tabularnewline
\end{tabular}\caption{\label{tab:LiH}The number of terms in the global group $N_{global}$
the number of non-empty groups $N_{group}$ and the maximum $N_{max}$,
minimum $N_{min}$ and average $N_{ave}$ number of terms in these
non-empty groups in $\phi_{1}$, $\phi_{2}$, $\mathcal{V}_{tot}\phi_{1}$
and $\mathcal{V}_{tot}\phi_{2}$ at the final step.}
\par\end{centering}
\end{table}

The resulting number of terms is again much smaller (about 100 times
smaller) than that needed using MADNESS software and by a factor $\approx2$
larger than that in \cite{BEY-HAU:2013}, where the number of terms
to represent the Slater-type orbitals is $1282$ and $1327$.

\begin{figure}[h]
\begin{centering}
\includegraphics[scale=0.55]{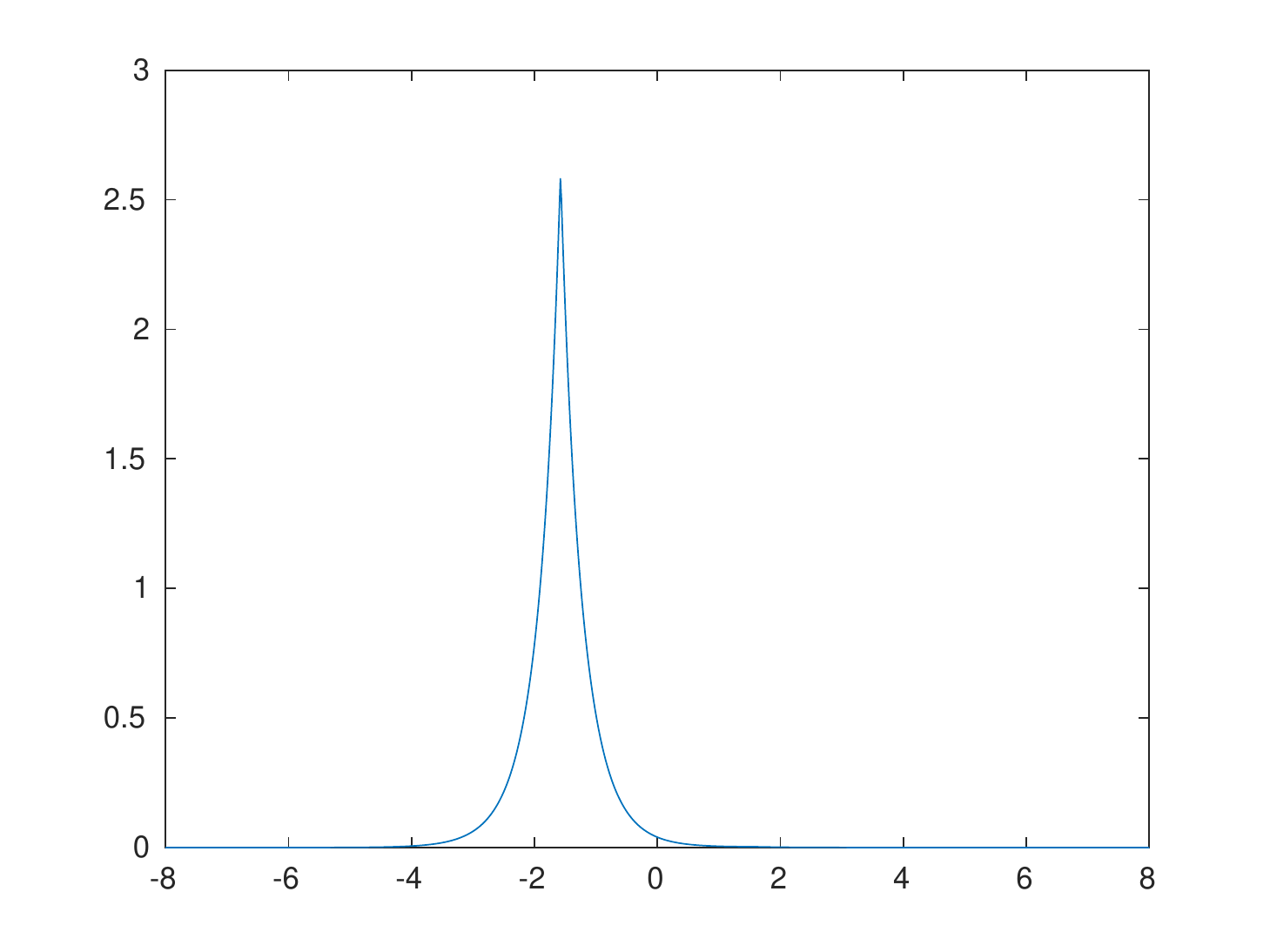}\includegraphics[scale=0.55]{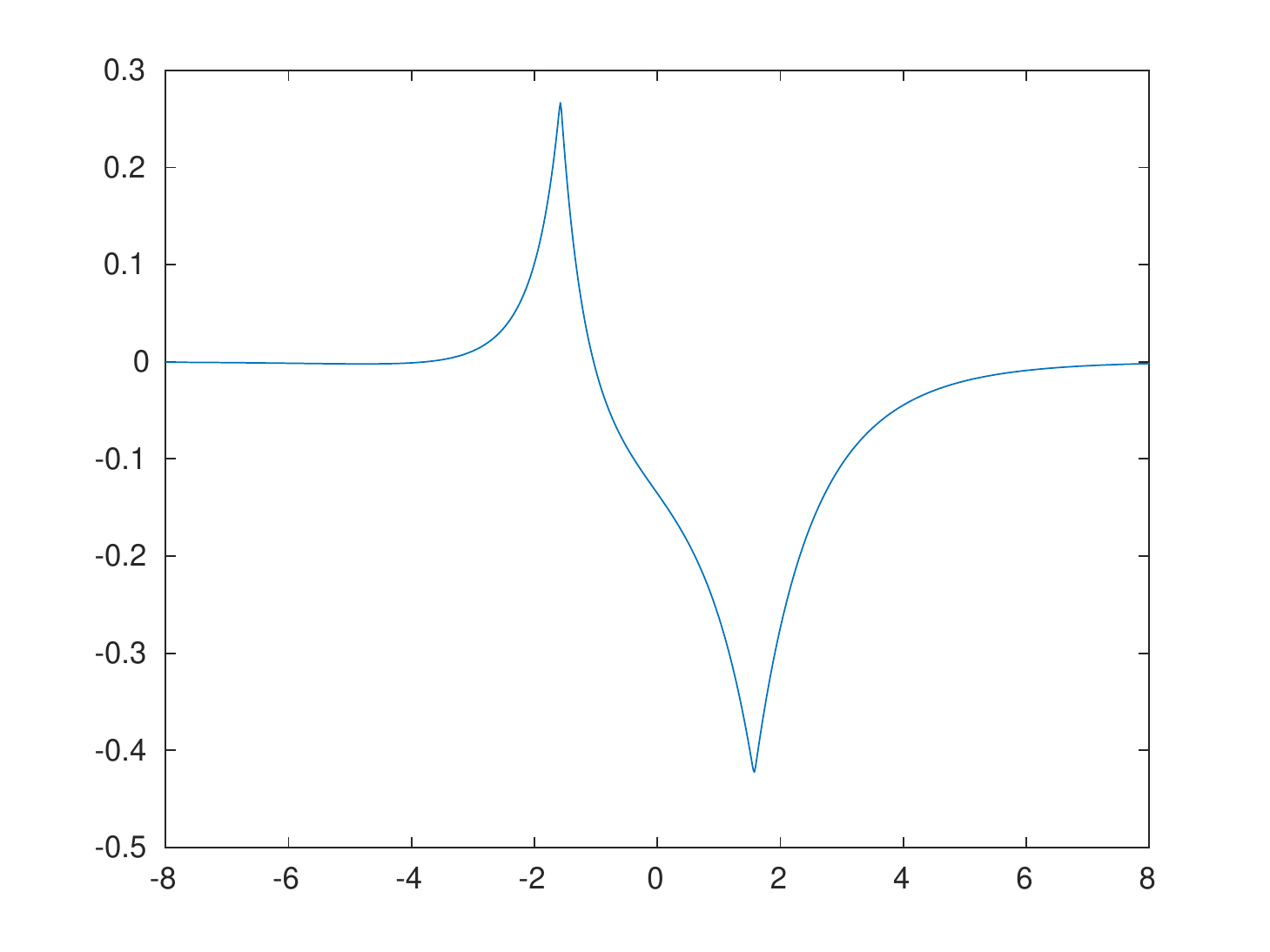}
\par\end{centering}
\caption{\label{fig:lithium}Plot of the spatial orbitals for lithium hydride,
$\mbox{LiH}$, $\phi_{1}\left(\mathbf{r}\right)$ (left) and $\phi_{2}\left(\mathbf{r}\right)$
(right) on the line $\mathbf{r}=\left(x,0,0\right)$ connecting the
two nucleus locations $\mathbf{R}_{1}=\left(-3.15/2,0,0\right)$ and
$\mathbf{R}_{2}=\left(3.15/2,0,0\right)$.}
\end{figure}

In Figure~\ref{fig:lithium} we display the spatial orbitals $\phi_{1}\left(\mathbf{r}\right)$
and $\phi_{2}\left(\mathbf{r}\right)$ on the line $\mathbf{r}=\left(x,0,0\right)$
connecting the two nuclei locations $\mathbf{R}_{1}$ and $\mathbf{R}_{2}$.

\begin{rem}
We implemented our code in Fortran90 and compiled it with the Intel
Fortran Compiler version 19.0.1.144. All computations were performed
on a single core (without any parallelization) of an Intel i7-6700
CPU at 3.4 GHz on a 64-bit Linux workstation with 64 GB of RAM. Currently
it takes about $119$ seconds to solve the Hartree\textendash Fock
equations for $\mbox{HeH}^{+}$ and $999$ seconds to solve the equations
for $\mbox{LiH}$. We made no attempt to optimize our implementation.
While the reduction algorithm that splits Gaussian mixtures into groups
(see Section~\ref{sec:Reduction-algorithm}) is trivially parallel,
we did not take advantage of this property in the current version
of our code. Since, for small molecules, the main cost is in reduction,
we expect a speedup factor of about $20$, given the number of groups
in Table~\ref{tab:LiH}. Such an acceleration factor should bring
the timing close to that of MADNESS. We plan to do a careful speed
comparison with existing methods separately.
\end{rem}

~
\begin{rem}
The reduction algorithm we used \cite[Algorithm 1]{BE-MO-YA:2019}
is capable of achieving near full double precision accuracy if certain
computations are carried out in quadruple precision. In addition to
solving (\ref{eq:helium eqn}) and (\ref{eq:lithium eqn}) using the
reduction algorithm implemented in double precision and obtaining
orbital energies $\approx10^{-6}$, we continued iterations and switched
to the reduction algorithm implemented in quadruple precision to improve
the accuracy. We observed that, within several additional iterations,
the orbital energies improved by gaining additional accurate digits.
However, as pointed out in \cite{BE-MO-YA:2019}, the quadruple precision
implementation is more than $10$ times slower than the double precision
version (if the same accuracy is required). We do not include these
results since, by modifying the reduction algorithm, it should be
possible to avoid using quadruple precision. However, we note that,
in most practical computations, an accuracy of $\approx10^{-6}$ is
typically sufficient. 
\end{rem}

\section{\label{sec:Conclusions-and-further-work}Conclusions and further
work}

We presented a new adaptive algorithm for computing the electronic
structure of molecules. Our approach is novel since we do not select
a basis in advance and let the iteration for solving a system of integral
equations identify the necessary functions to achieve a desired accuracy.
The reduction algorithm developed in \cite{BE-MO-YA:2019} selects
the ``best'' subset of linearly independent terms after each operation
that generates a large number of terms. These terms are Gaussians
(identified by their parameters) and, since potentials and Green's
functions are also represented via Gaussian mixtures, all integrals
are evaluated explicitly. Hence, in the algorithm, we only update
the parameters of Gaussian mixtures.

We introduce in this paper a hierarchical subdivision of terms of
Gaussian mixtures so that the computational cost of reduction is proportional
to the number of nuclear centers. Although we did not implement a
parallel version of our approach, it is clear that such an implementation
will have a significant impact on the speed since reduction within
each group of terms is independent of the other groups and computations
are also mostly independent between orbitals. 

We note that our method is less technical than that used in MADNESS
(see \cite{H-B-B-C-F-F-G-etc:2016}) and in \cite{BEY-HAU:2013} since
our only tool is the reduction algorithm. We plan further work on
our approach to verify its performance on large molecules.

\section*{Acknowledgement}

We thank the anonymous reviewers for their comments and encouragement.

\bibliographystyle{plain}

%%\bibliography{common}

\end{document}